\documentclass{amsart}%
\usepackage{amsthm,amstext,amscd,amsmath,amsfonts,amssymb}
\usepackage{amsmath}
\theoremstyle{plain}
\newtheorem{theorem}{Theorem}
\newtheorem{lemma}[theorem]{Lemma}

\theoremstyle{definition}
\newtheorem{definition}[theorem]{Definition}

\theoremstyle{remark}

\renewcommand{\Box}{\square}

\begin{document}

\title{Regular graded algebras and obstructed categories with duality}
\author{Steven Duplij}
\address{Kharkov National University, Svoboda Sq. 4, Kharkov 61077, Ukraine}
\email{Steven.A.Duplij@univer.kharkov.ua}
\urladdr{http://gluon.physik.uni-kl.de/\~{}duplij}
\author{Wladyslaw Marcinek}
\address{Institute of Theoretical Physics, University of Wroclaw,\\Pl.
Maxa Borna 9, 50-204 Wroc{l}aw, Poland}
\email{wmar@ift.univ.wroc.pl}
 \maketitle
\begin{abstract}
Regular and higher regular graded algebras (in simplest case satisfying Von
Neumann regularity $\Theta_{1}\Theta_{2}\Theta_{1}=\Theta_{1}$ instead of
anticommutativity) are introduced and their properties are studied. They are
described in terms of obstructed categories with nonclosed invertible and
noninvertible morphisms for which generalized (obstructed) functors and
natural transformations are given. Regular algebras and bialgebras are
considered with examples. Corresponding regularization of the cross product
and Wick algebras is made.

\end{abstract}

\section{Introduction}

The concept of higher regularization was introduced in the abstract way in the
supermanifold theory \cite{duplij,dup18} and then considered for general
morphisms \cite{dup/mar1,dup/mar}, which leaded to regularization of
categories and Yang-Baxter equation \cite{dup/mar2,dup/mar3,dup/mar5}. This
concept is close to the some generalizations of category theory
\cite{baez,bae/dol3,bae/dol1,gra/ozi}, it is connected with weak bialgebras
\cite{nik/vai,nil1}, squared Hopf algebras \cite{lyu2} and weak Hopf algebras
\cite{fangli3,dup/li1,dup/li2} and plays an important role in topological
quantum field theories \cite{cra/yet,bae/dol}. Here we introduce higher
regular graded algebras and study them in terms of obstructed categories,
which can clear and concrete understanding of the concept of semisupermanifold
\cite{duplij}.

Let $\mathcal{A}=\mathcal{A}_{\overline{0}}\oplus\mathcal{A}_{\overline{1}}$
be an arbitrary superalgebra equipped with an associative multiplication
$\mathrm{m}:\mathcal{A\otimes\mathcal{A}\rightarrow\mathcal{A}}$ [16]. We use
the notation $\mathrm{m}\left(  a\otimes b\right)  \equiv ab$ for the
multiplication (we call it internal product). Denote by $\left|  a\right|  $
the parity of $a\in\mathcal{A}$. Then we have the relation $\left|  ab\right|
=\left|  a\right|  +\left|  b\right|  ($mod$\,2)$ for the multiplication in
$\mathcal{A}$. The multiplication $\mathrm{m}$ is graded homogeneous mapping
\begin{equation}
\mathrm{m}_{\overline{ij}}:\mathcal{A}_{\overline{i}}\otimes\mathcal{A}%
_{\overline{j}}\rightarrow\mathcal{A}_{\overline{i}+\overline{j}%
(\text{mod}\,2)},
\end{equation}
where $\mathrm{m}_{\overline{ij}}$ is the restriction of $\mathrm{m}$ to
$\mathcal{A}_{\overline{i}}\otimes\mathcal{A}_{\overline{j}}$ and
$\overline{i},\overline{j}=0,1$. It is obvious that the mapping $\mathrm{m}%
_{\overline{00}}:\mathcal{A}_{\overline{0}}\otimes\mathcal{A}_{\overline{0}%
}\rightarrow\mathcal{A}_{\overline{0}}$ is a multiplication turning
$\mathcal{A}_{\overline{0}}$ into a closed associative algebra. Mappings
$\mathrm{m}_{\overline{01}}$ and $\mathrm{m}_{\overline{10}}$ define on
$\mathcal{A}_{\overline{1}}$ an $\mathcal{A}_{\overline{0}}$-module structure.
The mapping $\mathrm{m}_{\overline{11}}:\mathcal{A}_{\overline{1}}%
\otimes\mathcal{A}_{\overline{1}}\rightarrow\mathcal{A}_{\overline{0}}$ is an
$\mathcal{A}_{\overline{0}}$-valued bilinear form on $\mathcal{A}%
_{\overline{1}}$. This means that the product $ab$ is an even element of
$\mathcal{A}$, $\left|  ab\right|  =0\,($mod$\,2)$ for arbitrary two odd
elements $a,b\in\mathcal{A}_{\overline{1}}$. Supercommutativity \cite{berezin}
means that $ab\mathbf{-}\left(  -1\right)  ^{\left|  a\right|  \left|
b\right|  }ba\mathbf{=}0$. Let us consider a graded algebra $\mathcal{A}$
generated by a set $E\subset\mathcal{A}_{\overline{1}}$ of odd elements. In
this case $\mathcal{A}_{\overline{0}}$ is the image of $\mathrm{m}%
_{\overline{1}\overline{1}}$
\begin{equation}
\mathcal{A}_{\overline{0}}:=\mbox{Im} (\mathrm{m}_{\overline{1}\overline{1}}).
\end{equation}
Note that $\mathcal{A}_{\overline{1}}$ is also closed under the triple product
$\mathrm{m}_{\overline{1}\overline{1}\overline{1}}:\mathcal{A}_{\overline{1}%
}\otimes\mathcal{A}_{\overline{1}}\otimes\mathcal{A}_{\overline{1}}%
\rightarrow\mathcal{A}_{\overline{1}}$, where
\begin{equation}
\mathrm{m}_{\overline{1}\overline{1}\overline{1}}:=\mathrm{m}_{\overline
{0}\overline{1}}\circ(\mathrm{m}_{\overline{1}\overline{1}}\times
\operatorname*{id})=\mathrm{m}_{\overline{1}\overline{0}}\circ
(\operatorname*{id}\times\mathrm{m}_{\overline{1}\overline{1}}).
\end{equation}

\section{Regular graded algebras}

Here we consider graded algebras $\mathcal{A}^{reg}$ as a direct product of
superalgebras with new multiplication $\left(  \cdot\right)  $ (we call it
external product) between copies satisfying the following regularity
relations
\begin{align}
\mathbf{a}_{1}\mathbf{\cdot a}_{2}\mathbf{\cdot a}_{1}  &  =\mathbf{a}%
_{1},\;\;\;\;\mathbf{a}_{2}\mathbf{\cdot a}_{1}\mathbf{\cdot a}_{2}%
=\mathbf{a}_{2},\label{1}\\
\mathbf{a}_{1}  &  \in\mathcal{A}_{1},\;\;\mathbf{a}_{2}\in\mathcal{A}%
_{2},\;\;\mathcal{A}^{reg}=\mathcal{A}_{1}\times\mathcal{A}_{2},\nonumber
\end{align}
or their higher order generalizations \cite{dup/mar1,dup/mar}. We consider a
class of graded algebras as a free product of $n$ copies of a superalgebra
$\mathcal{A}$ modulo the regularity relations (\ref{1}) and call them
\textit{regular graded algebras}.

Recall (\cite{bor/mar} and refs. therein) that the free product of algebras
$\mathcal{A}_{1}$ and $\mathcal{A}_{2}$ is the algebra $\mathcal{A}_{1}%
\times\mathcal{A}_{2}$ formed by all formal finite sums of monomials of the
form $\mathbf{a}_{1}\cdot\mathbf{b}_{2}\cdot\mathbf{a}_{1}\ldots$ or
$\mathbf{b}_{1}\cdot\mathbf{a}_{2}\cdot\mathbf{b}_{1}\ldots$ , where
$\mathbf{a}_{i},\mathbf{b}_{i}\in\mathcal{A}_{i}$, $i=1,2$, without relations
between them.

Let us consider higher generalizations of (\ref{1}).

\begin{definition}
An element $\mathbf{a}_{1}\in\mathcal{A}_{1}$ for which there are elements
$\mathbf{a}_{i}\in\mathcal{A}_{i}$, where $i=2,\ldots n$ satisfying the
following condition
\begin{equation}%
\begin{array}
[c]{c}%
\mathbf{a}_{1}\cdot\mathbf{a}_{2}\cdot\ldots\cdot\mathbf{a}_{n}\cdot
\mathbf{a}_{1}=\mathbf{a}_{1},
\end{array}
\end{equation}
and its all cyclic permutations is said to be $n$-\textit{regular} in
$\mathcal{A}=\mathcal{A}_{1}\times\mathcal{A}_{2}\times\ldots\mathcal{A}_{n}$.
The set of all $n$-regular elements of $\mathcal{A}$ is denoted by
$\mathcal{A}^{reg}=Reg^{(n)}(\mathcal{A})$.
\end{definition}

The above $n$-regularity condition is essential for all noninvertible elements
of the algebra $\mathcal{A}$.

\begin{definition}
If the product of arbitrary two $n$-regular elements $\mathbf{a},\mathbf{b}\in
Reg^{(n)}(\mathcal{A})$ is also $n$-regular $\mathbf{a\cdot b}\in
Reg^{(n)}(\mathcal{A})$, then the multiplication $\mathrm{m}^{reg}%
(\mathbf{a}\otimes\mathbf{b})=\mathbf{a\cdot b}$ is said to be $n$-regular.
\end{definition}

The (super-)algebra $\mathcal{A}$ generating by $n$-regular elements and
equipped with $n$-regular multiplication is said to be $n$-regular.

\section{Regular odd graded algebras}

As a clear example we consider graded algebras $\mathcal{A}=\mathcal{A}%
_{\overline{0}}\oplus\mathcal{A}_{\overline{1}}$ generated by $n$-regular odd
elements $Reg(\mathcal{A}_{\overline{1}})$ which we call \textit{regular odd
graded algebras}. Obviously, that for these algebras every monomial formed by
an external product of an even number of generators is an even element of
$\mathcal{A}$, and the monomial with an odd number of generators is an odd
element. The superalgebra $\mathcal{A}$ contains all $\mathbb{K}$-linear
combinations of such monomials and external products.

We consider as an example an associative graded algebra $\mathcal{A}%
=\mathcal{A}^{reg}\left(  \Theta_{1},\Theta_{2}\right)  $ generated by two
noncommuting odd generators $\Theta_{1}$ and $\Theta_{2}$ satisfying the
following relations
\begin{equation}
\Theta_{1}^{2}=\Theta_{2}^{2}=0.
\end{equation}
and
\begin{equation}
\Theta_{1}\Theta_{2}\Theta_{1}=\Theta_{1},\;\;\;\Theta_{2}\Theta_{1}\Theta
_{2}=\Theta_{2}. \label{rel}%
\end{equation}

So the graded algebra $\mathcal{A}=\mathcal{A}^{reg}\left(  \Theta_{1}%
,\Theta_{2}\right)  $ is regular, but \textit{not} anticummutative. Let us
introduce
\[
E_{12}=\Theta_{1}\Theta_{2},\;\;\;\;E_{21}=\Theta_{2}\Theta_{1}.
\]

It is easy to observe that $E_{12}$ and $E_{21}$ are idempotents, since
\begin{align*}
E_{12}E_{12}  & =\left(  \Theta_{1}\Theta_{2}\Theta_{1}\right)  \Theta
_{2}=\Theta_{1}\Theta_{2}=E_{12},\\
E_{21}E_{21}  & =\left(  \Theta_{2}\Theta_{1}\Theta_{2}\right)  \Theta
_{1}=\Theta_{2}\Theta_{1}=E_{21}.
\end{align*}

Note that idempotents $E_{ij}$ are orthogonal $E_{12}E_{21}=E_{21}E_{12}=0$
and associative, and therefore they form an idempotent semigroup
\cite{petrich2,she/ovs}. From (\ref{rel}) it follows that $E_{ij}$ are
two-sided units for $\Theta_{i}$, viz. $E_{12}\Theta_{1}=\Theta_{1}$,
$\Theta_{1}E_{21}=\Theta_{1}$, $E_{21}\Theta_{2}=\Theta_{2}$, $\Theta
_{2}E_{12}=\Theta_{2}$. A generic element in $\mathcal{A}^{reg}(\Theta
_{1},\Theta_{2})$ is
\begin{equation}
\mathbf{a}=a_{0}+a_{1}\Theta_{1}+a_{2}\Theta_{2}+a_{12}\Theta_{1}\Theta
_{2}+a_{21}\Theta_{2}\Theta_{1}\label{a}%
\end{equation}
The multiplication in $\mathcal{A}^{reg}(\Theta_{1},\Theta_{2})$ is given by
\begin{equation}%
\begin{array}
[c]{l}%
\mathbf{ab}=a_{0}b_{0}\\
+(a_{0}b_{1}+a_{1}b_{0}+\underline{a_{1}b_{21}+a_{12}b_{1}})\Theta_{1}\\
+(a_{0}b_{2}+a_{2}b_{0}+\underline{a_{2}b_{12}+a_{21}b_{2}})\Theta_{2}\\
+(a_{0}b_{12}+a_{1}b_{2}+a_{12}b_{0}+\underline{a_{12}b_{12}})\Theta_{1}%
\Theta_{2}\\
+(a_{0}b_{21}+a_{2}b_{1}+a_{21}b_{0}+\underline{a_{21}b_{21}})\Theta_{2}%
\Theta_{1}.
\end{array}
\label{ab}%
\end{equation}

Note that the underlined terms are absent in anticommutative algebras
\cite{berezin}, and so they appear \textit{only} due to regularity
(\ref{rel}). That considerably changes the whole machinery of regular graded
algebras computations as compared with anticommutative case \cite{berezin}.
For instance, the equation $\mathbf{ab}=0$ has the solution $\mathbf{b}%
_{zdiv}=b_{0}\left(  1-\Theta_{1}-\Theta_{2}\right)  $, which means that
elements $\mathbf{b}_{zdiv}$ are zero divisors for all elements of the algebra
$\mathcal{A}^{reg}(\Theta_{1},\Theta_{2})$. For comparison, in $\Lambda\left(
\theta_{1},\theta_{2}\right)  $ (see e.g. \cite{berezin}) we have
$\mathbf{b}_{zdiv}^{\Lambda}=k\left(  a_{2}\theta_{1}+a_{1}\theta_{2}\right)
+r\theta_{1}\theta_{2}$, where $k$ and $r$ are any in a given field. Then,
instead of well known formula for inverse in $\Lambda\left(  \theta_{1}%
,\theta_{2}\right)  $ as ($a_{0}\neq0$)
\[
\left(  a_{0}+a_{1}\theta_{1}+a_{2}\theta_{2}+a_{12}\theta_{1}\theta
_{2}\right)  ^{-1}=a_{0}^{-1}-a_{0}^{-2}\left(  a_{1}\theta_{1}+a_{2}%
\theta_{2}+a_{12}\theta_{1}\theta_{2}\right)
\]
we have in $\mathcal{A}^{reg}(\Theta_{1},\Theta_{2})$%
\begin{align*}
\mathbf{a}^{-1}  & =a_{0}^{-1}-D^{-1}[a_{1}\Theta_{1}+a_{2}\Theta_{2}\\
& +\left(  a_{12}\left(  1+\dfrac{a_{21}}{a_{0}}\right)  -\dfrac{a_{1}a_{2}%
}{a_{0}}\right)  \Theta_{1}\Theta_{2}+\left(  a_{21}\left(  1+\dfrac{a_{12}%
}{a_{0}}\right)  -\dfrac{a_{1}a_{2}}{a_{0}}\right)  \Theta_{2}\Theta_{1}],
\end{align*}
where $D=\left(  a_{0}+a_{12}\right)  \left(  a_{0}+a_{21}\right)  -a_{1}%
a_{2}$.

\section{Regular graded algebra representations}

The graded algebra $\mathcal{A}=\mathcal{A}^{reg}(\Theta_{1},\Theta_{2})$ can
be presented by a free product of two one dimensional Grassmann algebras
$\Lambda\left(  \theta_{1}\right)  $ and $\Lambda\left(  \theta_{2}\right)  $
modulo regularity relations (\ref{1}) (denote them by $\;\overset
{reg}{\thicksim}$ ) as%
\[
\mathcal{A}^{reg}(\Theta_{1},\Theta_{2})\cong\Lambda(\theta_{1})\times
\Lambda(\theta_{2})\diagup\;\overset{reg}{\thicksim}.
\]

That means $\Theta_{1}\Theta_{2}\Theta_{1}=\Theta_{1}\Longleftrightarrow
\theta_{1}\cdot\theta_{2}\cdot\theta_{1}=\theta_{1}$, where the external
$\left(  \cdot\right)  $-product of elements from different copies of one
dimensional Grassmann algebra is defined in (\ref{1}).

Let us consider $n$ copies of one-dimensional Grassmann algebra $\Lambda
\left(  \theta_{1}\right)  $. The $i$-th copy we denote by $\Lambda(\theta
_{n})$. Let us define a graded algebra $\mathcal{A}^{reg}(\Theta_{1}$
$,\Theta_{2}$ $,\ldots,\Theta_{n})$ as a free product of $n$ copies of
one-dimensional Grassmann algebras subject to the following relation
\begin{equation}%
\begin{array}
[c]{c}%
\Theta_{1}\Theta_{2}\Theta_{3}\ldots\Theta_{n}\Theta_{1}=\Theta_{1}%
\end{array}
\end{equation}
and its all cyclic permutations. So the graded algebra $\mathcal{A}%
^{reg}(\Theta_{1}$ $,\Theta_{2}$ $,\ldots,\Theta_{n})$ is regular. Thus the
graded algebra $\mathcal{A}=\mathcal{A}^{reg}(\Theta_{1}$ $,\Theta_{2}$
$,\ldots,\Theta_{n})$ can also be presented by a free product of $n$ one
dimensional Grassmann algebras $\Lambda\left(  \theta_{1}\right)  $ modulo
regularity relations (\ref{1}) (denote them by $\;\overset{reg}{\thicksim}$ )
as%
\[
\mathcal{A}^{reg}(\Theta_{1},\Theta_{2},\ldots,\Theta_{n})\cong\Lambda
(\theta_{1})\times\Lambda(\theta_{2})\times\ldots\times\Lambda(\theta
_{n})\diagup\;\overset{reg}{\thicksim}.
\]

For the $n$-regular graded algebra $\mathcal{A}^{reg}(\Theta_{1}$ $,\Theta
_{2}$ $,\ldots,\Theta_{n})$ we introduce a collection of right and left
multiplication operators $\mathsf{R}_{\Theta_{i}}$ and $\mathsf{L}_{\Theta
_{i}}$ by the formulas
\begin{equation}%
\begin{array}
[c]{c}%
\mathsf{R}_{\Theta_{i}}(x):=x\Theta_{i},\;\;\;\mathsf{L}_{\Theta_{i}%
}(x):=\Theta_{i}x,
\end{array}
\end{equation}
for every $x\in\mathcal{A}^{reg}(\Theta_{1}$ $,\Theta_{2}$ $,\ldots,\Theta
_{n})$ and $i=1,\ldots n$.

\begin{lemma}
\label{lemm-repr}The correspondence defines a representation of $\mathcal{A}%
^{reg}(\Theta_{1}$ $,\Theta_{2}$ $,\ldots,\Theta_{n})$ into itself
\begin{equation}%
\begin{array}
[c]{c}%
\mathcal{L}:\Theta_{i}\rightarrow\mathsf{L}_{\Theta_{i}}\in\operatorname*{End}%
(\mathcal{A}^{reg}(\Theta_{1},\Theta_{2},\ldots,\Theta_{n}))
\end{array}
\end{equation}
and%
\begin{equation}%
\begin{array}
[c]{c}%
\mathcal{R}:\Theta_{i}\rightarrow\mathsf{R}_{\Theta_{i}}\in\operatorname*{End}%
(\mathcal{A}^{reg}(\Theta_{1},\Theta_{2},\ldots,\Theta_{n}))
\end{array}
\end{equation}
defines an antirepresentation of $\mathcal{A}^{reg}(\Theta_{1}$ $,\Theta_{2}$
$,\ldots,\Theta_{n})$ into itself.
\end{lemma}

\textbf{Proof:} We calculate
\[%
\begin{array}
[c]{l}%
\mathsf{R}_{\Theta_{i}}\mathsf{R}_{\Theta_{j}}=\mathsf{R}_{\Theta_{i}%
}(\mathsf{R}_{\Theta_{j}}(x))=\mathsf{R}_{\Theta_{i}}\left(  x\Theta
_{j}\right)  =x\Theta_{j}\Theta_{i}=\mathsf{R}_{\Theta_{j}\Theta_{i}},\\
\mathsf{L}_{\Theta_{i}}\mathsf{L}_{\Theta_{j}}=\mathsf{L}_{\Theta_{i}%
}(\mathsf{L}_{\Theta_{j}}(x))=\mathsf{L}_{\Theta_{i}}\left(  \Theta
_{j}x\right)  =\Theta_{i}\Theta_{j}x=\mathsf{L}_{\Theta_{i}\Theta_{j}}.
\end{array}
\]

\section{Regular graded algebra decomposition}

Now we use Lemma \ref{lemm-repr} to construct the graded algebra decomposition
into linear spaces. First we consider the regular graded algebra
$\mathcal{A}^{reg}(\Theta_{1},\Theta_{2})$ and introduce 2 two--dimensional
linear spaces $\mathtt{X}_{1}$ and $\mathtt{X}_{2}$ over a field $\mathbb{K}$
\begin{equation}%
\begin{array}
[c]{c}%
\mathtt{X}_{1}:=lin_{\mathbb{K}}\{\Theta_{1},\Theta_{1}\Theta_{2}%
\}=lin_{\mathbb{K}}\{\Theta_{1},E_{12}\},\\
\mathtt{X}_{2}:=lin_{\mathbb{K}}\{\Theta_{2},\Theta_{2}\Theta_{1}%
\}=lin_{\mathbb{K}}\{\Theta_{2},E_{21}\}.
\end{array}
\end{equation}
Define two linear mappings $f_{1}:\mathtt{X}_{1}\rightarrow\mathtt{X}_{2}^{T}$
and $f_{2}:\mathtt{X}_{2}\rightarrow\mathtt{X}_{1}^{T}$ (where $T$ is
coordinate transposition, $T^{2}=1$) as a left multiplication by $\Theta_{2}$
and $\Theta_{1}$ respectively (because $\mathcal{L}$ is a representation). We
obtain
\begin{equation}%
\begin{array}
[c]{ll}%
f_{1}(\Theta_{1})=\mathsf{L}_{\Theta_{2}}\left(  \Theta_{1}\right)
=\Theta_{2}\Theta_{1}, & f_{1}(\Theta_{1}\Theta_{2})=\mathsf{L}_{\Theta_{2}%
}\left(  \Theta_{1}\Theta_{2}\right)  =\Theta_{2}\Theta_{1}\Theta_{2}%
=\Theta_{2},\\
f_{2}(\Theta_{2})=\mathsf{L}_{\Theta_{1}}\left(  \Theta_{2}\right)
=\Theta_{1}\Theta_{2}, & f_{2}(\Theta_{2}\Theta_{1})=\mathsf{L}_{\Theta_{1}%
}\left(  \Theta_{2}\Theta_{1}\right)  =\Theta_{1}\Theta_{2}\Theta_{1}%
=\Theta_{1},
\end{array}
\end{equation}
where regularity plays an essential role. Obviously we have
\begin{equation}%
\begin{array}
[c]{c}%
f_{1}\circ f_{2}\circ f_{1}=f_{1},\quad f_{2}\circ f_{1}\circ f_{2}=f_{2}.
\end{array}
\label{fff}%
\end{equation}
Define two mappings
\begin{equation}%
\begin{array}
[c]{c}%
e_{\mathtt{X}_{1}}:=f_{2}\circ f_{1},\quad e_{\mathtt{X}_{2}}:=f_{1}\circ
f_{2}.
\end{array}
\label{eff}%
\end{equation}
It is obvious that we have
\begin{equation}%
\begin{array}
[c]{c}%
e_{\mathtt{X}_{1}}\circ e_{\mathtt{X}_{1}}=e_{\mathtt{X}_{1}},\quad
e_{\mathtt{X}_{2}}\circ e_{\mathtt{X}_{2}}=e_{\mathtt{X}_{2}}.
\end{array}
\end{equation}

Observe that mappings $e_{\mathtt{X}_{i}}$ are identities only for regular
graded algebras $\mathcal{A}^{reg}(\Theta_{1},\Theta_{2})$. Thus we can
construct a category $\mathcal{C}_{\mathcal{A}}^{\left(  2\right)  }$ whose
objects are $\mathcal{C}_{\operatorname*{Ob}}:=\{\mathbb{K},\mathtt{X}%
_{1},\mathtt{X}_{2},\;\mathtt{X}_{1}\oplus\mathtt{X}_{2}\}$ and whose
morphisms are given as all compositions of mappings $\{f_{1},f_{2}%
,e_{\mathtt{X}_{1}},e_{\mathtt{X}_{2}}\}$ and identity morphisms.

For an $n$-regular graded algebra $\mathcal{A}=\mathcal{A}^{reg}(\Theta_{1}$
$,\Theta_{2}$ $,\ldots,\Theta_{n})$ generated by a set of regular odd
generators $\{\Theta_{i};i=1,\ldots,n\}$ we introduce the following
decomposition
\begin{equation}%
\begin{array}
[c]{c}%
\mathcal{A}=\bigoplus_{i=1}^{n}\mathtt{X}_{i},
\end{array}
\end{equation}
such that $\mathsf{L}_{\Theta_{i}}(x)\in\mathtt{X}_{i+1}$ for every
$x\in\mathtt{X}_{i}$, and $i=1,\ldots,n-1$. The space $\mathtt{X}_{1}$ is the
$\mathbb{K}$-linear span of all monomials in generators ended by $\Theta_{i}$.
For $i=2,3,\ldots,n$ we define
\begin{equation}%
\begin{array}
[c]{c}%
\mathtt{X}_{i}:=lin_{\mathbb{K}}\{\mathsf{L}_{\Theta_{i-1}}\circ\cdots
\circ\mathsf{L}_{\Theta_{1}}(x):x\in\mathtt{X}_{1};i=2,\ldots,n\}.
\end{array}
\end{equation}
and obtain
\begin{equation}%
\begin{array}
[c]{l}%
\mathtt{X}_{2}:=lin_{\mathbb{K}}\{\mathsf{L}_{\Theta_{1}}(x):x\in
\mathtt{X}_{1}\},\\
\mathtt{X}_{3}:=lin_{\mathbb{K}}\{\mathsf{L}_{\Theta_{2}}\circ\mathsf{L}%
_{\Theta_{1}}(x):x\in\mathtt{X}_{1}\},\\
\hspace{2cm}\vdots\hspace{3cm}\vdots\\
\mathtt{X}_{n}:=lin_{\mathbb{K}}\{\mathsf{L}_{\Theta_{n-1}}\circ\cdots
\circ\mathsf{L}_{\Theta_{1}}(x):x\in\mathtt{X}_{1}\}.
\end{array}
\end{equation}

We use operators $\mathsf{L}_{a}$ in order to introduce a collection of linear
mappings $(\mathtt{X},f_{i}):=\{f_{i}:\mathtt{X}_{i}\rightarrow\mathtt{X}%
_{i+1};i=1,\ldots,n\}$. We define these mappings by the following formula
\begin{equation}%
\begin{array}
[c]{c}%
f_{i}(x):=\mathsf{L}_{\Theta_{i}}(x)=\Theta_{i}x
\end{array}
\end{equation}
for every $x\in\mathtt{X}_{i}$, $i=1,\ldots,n-1$. One can calculate that we
have the higher regularity relation
\begin{equation}%
\begin{array}
[c]{c}%
f_{1}\circ f_{2}\circ\cdots\circ f_{n}\circ f_{1}=f_{1}%
\end{array}
\label{fn}%
\end{equation}
and corresponding cyclic permutations. In this case there is also a specific
category $\mathcal{C}_{\mathcal{A}}^{\left(  n\right)  }$ which contains all
spaces and mappings considered above.

\section{Categories and obstructions}

\subsection{Algebras via categories}

Recall (see e.g. \cite{maclane1}) that a \textit{category} $\mathcal{C}%
=(\mathcal{C}_{\operatorname*{Ob}},\mathcal{C}_{\operatorname*{Mor}%
},\mathsf{c})$ contains

(i) a collection $\mathcal{C}_{\operatorname*{Ob}}$ of objects;

(ii) a collection $\mathcal{C}_{\operatorname*{Mor}}$ of morphisms (arrows)%

\begin{equation}
\mathcal{C}_{\operatorname*{Mor}}=\bigcup_{X,Y\in\mathcal{C}%
_{\operatorname*{Ob}}}\mathcal{C}(X,Y)
\end{equation}

(iii) an associative composition $\mathsf{c}$ of morphisms%

\begin{equation}
\mathsf{c}:\mathcal{C}(X,Y)\times\mathcal{C}(Y,\mathcal{W})\rightarrow
\mathcal{C}(X,\mathcal{W}) \label{com}%
\end{equation}

The collection $\mathcal{C}_{\operatorname*{Mor}}$ is the union of mutually
disjoint sets $\mathcal{C}(X,Y)$ of morphisms $X\rightarrow Y$ from $X$ to $Y$
defined for every pair of objects $X,Y\in\mathcal{C}_{\operatorname*{Ob}}$. It
may happen that for a pair $X,Y\in\mathcal{C}_{\operatorname*{Ob}}$ the set
$\mathcal{C}(X,Y)$ is empty.

An \textit{opposite} (or \textit{dual}) category of a category $\mathcal{C}%
=(\mathcal{C}_{\operatorname*{Ob}},\mathcal{C}_{\operatorname*{Mor}%
},\mathsf{c})$ is a category $\mathcal{C}^{op}=(\mathcal{C}%
_{\operatorname*{Ob}},\mathcal{C}_{\operatorname*{Mor}}^{op},\mathsf{c}^{op})$
equipped with the same collection of objects $\mathcal{C}_{\operatorname*{Ob}%
}$ as the category $\mathcal{C}$ but with reversed all arrows
\begin{equation}
\mathcal{C}^{op}(X,Y)\equiv\mathcal{C}(Y,X).
\end{equation}
and reversed composition $\mathsf{c}^{op}$. If $\mathfrak{D}_{\mathcal{C}}$ is
a diagram built from objects and morphisms of the category $\mathcal{C}$, then
the same diagram with reversed all arrows $\mathfrak{D}_{\mathcal{C}^{op}}$ is
said to be \textit{dual} to $\mathfrak{D}_{\mathcal{C}}$.

Let $\mathcal{C}$ and $\mathcal{D}$ be two categories. A \textit{functor}
$\mathfrak{F}:\mathcal{C}\rightarrow\mathcal{D}$ of $\mathcal{C}$ into
$\mathcal{D}$ is a pair of maps
\begin{equation}%
\begin{array}
[c]{cc}%
\mathfrak{F}_{\operatorname*{Ob}}:\mathcal{C}_{\operatorname*{Ob}}%
\rightarrow\mathcal{D}_{\operatorname*{Ob}}, & \mathfrak{F}%
_{\operatorname*{Mor}}:\mathcal{C}_{\operatorname*{Mor}}\rightarrow
\mathcal{D}_{\operatorname*{Mor}}%
\end{array}
\end{equation}
such that
\begin{equation}
\mathfrak{F}_{\operatorname*{Mor}}(\varphi\circ\psi)=\mathfrak{F}%
_{\operatorname*{Mor}}(\varphi)\circ\mathfrak{F}_{\operatorname*{Mor}}(\psi)
\end{equation}
for every morphisms $\varphi:Y\longrightarrow Z$ and $\psi:X\longrightarrow Y$
of $\mathcal{C}$. The generalization to multifunctors is obvious. For instance
an $m$--ary functor $\mathfrak{F}:\mathcal{C}^{\times m}\longrightarrow
\mathcal{D}$ sends an $m$-tuple of objects of $\mathcal{C}$ into an object of
$\mathcal{D}$.

A \textit{natural transformation} of functors $\mathfrak{s}:\mathfrak{F}%
\rightarrow\mathfrak{G}$ of $\mathfrak{F}$ into $\mathfrak{G}$ is a collection
of morphisms
\begin{equation}
\mathfrak{s}=\{\mathfrak{s}_{X}:\mathfrak{F}(X)\rightarrow\mathfrak{G}%
(X),\;\;X\in\mathcal{C}_{\operatorname*{Ob}}\},
\end{equation}
such that
\begin{equation}
\mathfrak{s}_{Y}\circ\mathfrak{F}(\psi)=\mathfrak{G}(\psi)\circ\mathfrak{s}%
_{X}%
\end{equation}
for every morphism $\psi:X\longrightarrow Y$ of $\mathcal{C}$. The set of all
natural transformations of $\mathfrak{F}$ into $\mathfrak{G}$ we denote by
$\mathfrak{Nat}(\mathfrak{F},\mathfrak{G})$. If $\mathfrak{F}\equiv
\mathfrak{G}$, then we say that the natural transformation $\mathfrak{s}%
:\mathfrak{F}\longrightarrow\mathfrak{G}$ is a natural transformation of
$\mathfrak{F}$ into itself.

We can use functors and natural transformations in order to describe certain
algebraic structures in categories. Let $\mathcal{C}$ be a category and let
$\mathfrak{F}$ be a functor $\mathfrak{F}:\mathcal{C}\rightarrow\mathcal{C}$.
A binary multiplication is a natural transformation $\mathrm{m}:\mathfrak{F}%
\times\mathfrak{F}\rightarrow\mathfrak{F}$ satisfying the associativity
condition $\mathrm{m}\circ(\mathrm{m}\times\mathfrak{F})=\mathrm{m}%
\circ(\mathfrak{F}\times\mathrm{m})$. Similarly, a comultiplication is a
natural transformation $\triangle:\mathfrak{F}\rightarrow\mathfrak{F}%
\times\mathfrak{F}$ satisfying the corresponding coassociativity condition
$(\triangle\times\mathfrak{F})\circ\triangle=(\mathfrak{F}\times
\triangle)\circ\triangle$. A monoidal category $\mathcal{C}\equiv
\mathcal{C}(\otimes,\mathbb{K})$ is a category $\mathcal{C}$ equipped with a
monoidal operation (a bifunctor) $\otimes:\mathcal{C}\times\mathcal{C}%
\rightarrow\mathcal{C}$, a unit object $\mathbb{K}$ satisfying some known
axioms \cite{joy/str,yet}.

We define a (left) $\vee$--operation
\begin{equation}
X^{\vee}:=\mathcal{C}(X,\mathbb{K})
\end{equation}
such that we have
\begin{equation}
X^{\vee\vee}=X,\quad(X\otimes Y)^{\vee}=Y^{\vee}\otimes X^{\vee}%
\end{equation}
for every $X,Y\in\mathcal{C}_{\operatorname*{Ob}}$. If $\psi:X\rightarrow Y$
is a morphism in $\mathcal{C}$, then there is also a corresponding
\textit{dual morphism} $\psi^{\vee}:Y^{\vee}\rightarrow X^{\vee}$
\begin{equation}
\psi^{\vee}(v^{\vee}):=u^{\vee},\;u\in X,v\in Y\;\mbox{if and only if}
\;\psi(u)=v,\;u^{\vee}\in X^{\vee},v^{\vee}\in Y^{\vee}.
\end{equation}

The corresponding category $\mathcal{C}$ is called a monoidal category with
duality \cite{joy/str,yet}.

A class $\mathsf{g}(\mathcal{C}):=\{\mathsf{g}_{X}\equiv\langle\ldots
|\ldots\rangle_{X}:X\in\mathcal{C}_{\operatorname*{Ob}}\}$ of $\mathbb{K}%
$--valued mappings $\mathsf{g}_{X}:X^{\vee}\otimes X\rightarrow\mathbb{K}$
together with a class $\mathsf{h}(\mathcal{C}):=\{\mathsf{h}_{X}%
:\mathbb{K}\rightarrow X\otimes X^{\vee}\}$ such that the following diagrams
\[%
\begin{array}
[c]{ccc}
& \operatorname*{id}\nolimits_{Y^{\vee}}\otimes\mathsf{g}_{X}\otimes
\operatorname*{id}\nolimits_{Y} & \\
Y^{\vee}\otimes X^{\vee}\otimes X\otimes Y & \rightarrow & Y^{\vee}\otimes Y\\
&  & \\
\parallel &  & \downarrow\mathsf{g}_{Y}\\
&  & \\
(X\otimes Y)^{\vee}\otimes(X\otimes Y) & \rightarrow & \mathbb{K}\\
& \mathsf{g}_{X\otimes Y} &
\end{array}
\]

and%

\[%
\begin{array}
[c]{rcl}
& \operatorname*{id}{}_{Y^{\vee}}\otimes\psi\otimes\operatorname*{id}%
{}_{X^{\vee}} & \\
Y^{\vee}\otimes X\otimes X^{\vee} & \rightarrow & Y^{\vee}\otimes Y\otimes
X^{\vee}\\
&  & \\
\operatorname*{id}{}_{Y^{\vee}}\otimes\mathsf{h}_{X}\uparrow &  &
\downarrow\mathsf{g}_{Y}\otimes\operatorname*{id}{}_{X^{\vee}}\\
&  & \\
Y^{\vee} & \rightarrow & X^{\vee}\\
& \psi^{\vee} &
\end{array}
\]
is said to be a pairing in $\mathcal{C}$ \cite{yet}.

\subsection{Obstructed categories}

Let $\mathcal{C}$ be a category with invertible and noninvertible morphisms
\cite{dup/mar1,dup/mar} and an equivalence. By an equivalence in $\mathcal{C}$
we mean a class of morphisms $\mathcal{C}^{\mathrm{P}}=\bigcup_{X,Y\in
\mathcal{C}_{\operatorname*{Ob}}}\mathcal{C}^{\mathrm{P}}(X,Y)$, where
$\mathcal{C}^{\mathrm{P}}(X,Y)$ is a subset of invertible morphisms in
$\mathcal{C}(X,Y)$ satisfying the property $\mathrm{P}$. Two objects $X,Y$ of
the category $\mathcal{C}$ are equivalent, if and only if there is an morphism
$X\overset{\psi}{\longrightarrow}Y$ in $\mathcal{C}^{\mathrm{P}}(X,Y)$ such
that $\psi^{-1}\circ\psi=\operatorname*{id}{}_{X}$ and $\psi\circ\psi
^{-1}=\operatorname*{id}\nolimits_{Y}$.

\begin{definition}
A sequence of (invertible and noninvertible) morphisms
\begin{equation}%
\begin{array}
[c]{c}%
X_{1}\overset{\psi_{1}}{\longrightarrow}X_{2}\overset{\psi_{2}}%
{\longrightarrow}\cdots\overset{\psi_{n-1}}{\longrightarrow}X_{n}\overset
{\psi_{n}}{\longrightarrow}X_{1}%
\end{array}
\label{regu}%
\end{equation}
such that there is an (endo-)morphism $\mathbf{e}_{X_{1}}^{\left(  n\right)
}:X_{1}\longrightarrow X_{1}$ defined by
\begin{equation}%
\begin{array}
[c]{c}%
\mathbf{e}_{X_{1}}^{\left(  n\right)  }:=\psi_{n}\circ\cdots\circ\psi_{2}%
\circ\psi_{1}%
\end{array}
\label{egu}%
\end{equation}
and subjects to the following relation
\begin{equation}
\psi_{1}\circ\psi_{n}\circ\cdots\circ\psi_{2}\circ\psi_{1}=\psi_{1}%
\end{equation}
is said to be a \textit{regular }$n$\textit{-cocycle} on $\mathcal{C}$ and it
is denoted by $(X,\psi)$.
\end{definition}

Two mappings $f_{1}$ and $f_{2}$ defined in (\ref{fff}) form a regular
$2$-cocycle on $\mathcal{C}_{\mathcal{A}}^{\left(  2\right)  }$. Mappings
$f_{1},\ldots,f_{n}$ from (\ref{fn}) form a regular $n$-cocycle on
$\mathcal{C}_{\mathcal{A}}^{\left(  n\right)  }$.

\begin{definition}
The morphism $\mathbf{e}_{X_{1}}^{(n)}$ is said to be an \textit{obstruction}
of the object $X_{1}$ corresponding to the regular $n$-cocycle $\psi=(\psi
_{1},\ldots\psi_{n})$ on $\mathcal{C}$. The obstruction (endo-)morphisms
$\mathbf{e}_{X_{i}}^{\left(  n\right)  }:X_{i}\longrightarrow X_{i}$
corresponding for $i=2,\ldots,n$ are defined by a suitable cyclic permutation
of above sequence.
\end{definition}

\begin{definition}
If the obstruction (endo-)morphisms $\mathbf{e}_{X}^{\left(  n\right)
}:X\longrightarrow X$ are defined for every object $X\in\mathcal{C}%
_{\operatorname*{Ob}}$, then the mapping $\varepsilon^{(n)}:X\rightarrow
\mathbf{e}_{X}^{(n)}\in\mathcal{C}(X,X)$ is called a regular $n$%
-\textit{cocycle obstruction structure} on $\mathcal{C}$.
\end{definition}

It is obvious that for an usual category all $\mathbf{e}_{X}^{\left(
n\right)  }$ are equal to identity $\mathbf{e}_{X}^{\left(  n\right)  }%
=Id_{X}$. We are interested with categories for which the obstruction
$\mathbf{e}_{X}^{\left(  n\right)  }$ differs from the identity (see e.g.
\cite{duplij}).

\begin{definition}
A category $\mathcal{C}$ equipped with a regular $n$-cocycle obstruction
structure $\varepsilon^{(n)}:X\in\mathcal{C}_{\operatorname*{Ob}}\rightarrow
e_{X}^{(n)}\in\mathcal{C}(X,X)$ such that $\mathbf{e}_{X}^{(n)}\neq
\operatorname*{id}\nolimits_{X}$ for some $X\in\mathcal{C}_{\operatorname*{Ob}%
}$ is called an \textit{obstructed} category. The minimum number $n=n_{obstr}$
for which it occurs will define a quantitative measure of obstruction
$n_{obstr}$.
\end{definition}

The categories $\mathcal{C}_{\mathcal{A}}^{\left(  2\right)  }$ and
$\mathcal{C}_{\mathcal{A}}^{\left(  n\right)  }$ considered above are obstructed.

Let $(Y,\varphi)$ a regular $n$-cocycle, i.e. a sequence of morphisms
$Y_{1}\overset{\varphi_{1}}{\longrightarrow}Y_{2}\overset{\varphi_{2}%
}{\rightarrow}\cdots\overset{\varphi_{n-1}}{\rightarrow}Y_{n}\overset
{\varphi_{n}}{\rightarrow}Y_{1}$ such that $\mathbf{e}_{Y_{1}}^{\left(
n\right)  }:=\varphi_{n}\circ\cdots\circ\varphi_{2}\circ\varphi_{1}$.

\begin{definition}
A sequence of morphisms $\alpha:=(\alpha_{1},\ldots,\alpha_{n})$ such that the
diagram
\begin{equation}%
\begin{array}
[c]{rccccccccl}
& X_{1} & \overset{\psi_{1}}{\longrightarrow} & X_{2} & \overset{\psi_{2}%
}{\longrightarrow}\cdots\overset{\psi_{n-1}}{\longrightarrow} & X_{n} &
\overset{\psi_{n}}{\longrightarrow} & X_{1} &  & \\
& \downarrow\alpha_{1} &  & \downarrow\alpha_{2} &  & \downarrow\alpha_{n} &
& \downarrow\alpha_{1} &  & \\
& Y_{1} & \overset{\varphi_{1}}{\longrightarrow} & Y_{2} & \overset
{\varphi_{2}}{\longrightarrow}\cdots\overset{\varphi_{n-1}}{\longrightarrow} &
Y_{n} & \overset{\varphi_{n}}{\longrightarrow} & Y_{1} &  &
\end{array}
\end{equation}
is commutative, is said to be a \textit{regular }$\mathit{n}$\textit{-cocycle
morphism} from $(X,\psi)$ to $(Y,\varphi)$ and is denoted by $\alpha
:(X,\psi)\rightarrow(Y,\varphi)$.
\end{definition}

Observe that for a regular $n$-cocycle morphism $\alpha:(X,\psi)\rightarrow
(Y,\varphi)$ we have the relation
\begin{equation}
\alpha_{1}\circ\mathbf{e}_{X_{1}}^{\left(  n\right)  }=\mathbf{e}_{Y_{1}%
}^{(n)}\circ\alpha_{1}.
\end{equation}

\begin{definition}
A regular $n$-cocycle obstruction morphism $\beta:(X,\psi)\rightarrow
(X^{\prime},\varphi)$ which sends the object $X_{i}$ into equivalent object
$X_{i}^{\prime}$ and morphism $\psi_{i}$ into $\varphi_{i}$ is said to be
\textit{obstruction }$n$\textit{-cocycle equivalence}. The corresponding
obstructions $\mathbf{e}_{X}^{(n)}$ and $\mathbf{e}_{X^{\prime}}^{(n)}$ are
also said to be \textit{equivalent}.
\end{definition}

Let $\mathcal{C}$ and $\mathcal{D}$ be two obstructed categories. The
morphisms $\mathbf{e}_{X}^{\left(  n\right)  }$ can be used to extend the
notion of functors. Let $\mathfrak{F}:\mathcal{C}\rightarrow\mathcal{D}$ be a
functor defined as usual as a pair of mappings $(\mathfrak{F}%
_{\operatorname*{Ob}},\mathfrak{F}_{\operatorname*{Mor}})$.

\begin{definition}
A set of category mappings $\mathfrak{F}^{\left(  n\right)  }:\mathcal{C}%
\rightarrow\mathcal{D}$ defined as $(\mathfrak{F}_{\operatorname*{Ob}}%
^{(n)},\mathfrak{F}_{\operatorname*{Mor}}^{(n)})$ such that
\begin{equation}%
\begin{array}
[c]{ccc}%
\mathfrak{F}_{\operatorname*{Ob}}^{(n)}\equiv\mathfrak{F}_{\operatorname*{Ob}%
}, & \mathfrak{F}_{\operatorname*{Mor}}^{\left(  n\right)  }\left(
\mathbf{e}_{\mathcal{C}_{\operatorname*{Ob}}}^{\left(  n\right)  }\right)
=\mathbf{e}_{\mathcal{D}_{\operatorname*{Ob}}(X)}^{(n)}, & \mathfrak{F}%
_{\operatorname*{Mor}}^{\left(  1\right)  }\equiv\mathfrak{F}%
_{\operatorname*{Mor}},
\end{array}
\label{fe}%
\end{equation}
where $\mathfrak{F}=\left(  \mathfrak{F}_{\operatorname*{Ob}},\mathfrak{F}%
_{\operatorname*{Mor}}\right)  $ is a functor, is said to be $n$%
-\textit{regular obstructed functor}.
\end{definition}

We see that all the standard definitions of a functor are not changed, but the
preservation of identity $\mathfrak{F}_{\operatorname*{Mor}}\left(
\operatorname*{id}\nolimits_{X}\right)  =\operatorname*{id}\nolimits_{Y}$,
where $Y=\mathfrak{F}_{\operatorname*{Ob}}\left(  X\right)  $, $X\in
\mathcal{C}_{\operatorname*{Ob}}$, $Y\in\mathcal{D}_{\operatorname*{Ob}}$ is
replaced by requirement of preservation of the morphisms $\mathbf{e}%
_{X}^{\left(  n\right)  }$,and so the obstructed functor $\mathfrak{F}%
^{\left(  n\right)  }$ becomes $n$-dependent.

\begin{lemma}
Let the sequence
\begin{equation}%
\begin{array}
[c]{c}%
X_{1}\overset{\psi_{1}}{\longrightarrow}X_{2}\overset{\psi_{2}}%
{\longrightarrow}\cdots\overset{\psi_{n-1}}{\longrightarrow}X_{n}\overset
{\psi_{n}}{\longrightarrow}X_{1}%
\end{array}
\label{regul}%
\end{equation}
be a regular $n$-cocycle in the category $\mathcal{C}$. If $\mathfrak{F}%
^{\left(  n\right)  }:\mathcal{C}\rightarrow\mathcal{D}$ is an $n$-regular
functor, then
\begin{equation}%
\begin{array}
[c]{c}%
\mathfrak{F}_{\operatorname*{Ob}}^{(n)}(X_{1})\overset{\mathfrak{F}%
_{\operatorname*{Mor}}^{(n)}(\psi_{1})}{\longrightarrow}\mathfrak{F}%
_{\operatorname*{Ob}}^{(n)}(X_{2})\overset{\mathfrak{F}_{\operatorname*{Mor}%
}^{(n)}(\psi_{2})}{\longrightarrow}\cdots\overset{\mathfrak{F}%
_{\operatorname*{Mor}}^{(n)}(\psi_{n-1})}{\longrightarrow}\mathfrak{F}%
_{\operatorname*{Ob}}^{\left(  n\right)  }(X_{n})\overset{\mathfrak{F}%
_{\operatorname*{Mor}}^{\left(  n\right)  }(\psi_{n})}{\longrightarrow
}\mathfrak{F}_{\operatorname*{Ob}}^{\left(  n\right)  }(X_{1})
\end{array}
\label{regula}%
\end{equation}
is a regular $n$-cocycle in the category $\mathcal{D}$, i.e. $Y_{1}%
\overset{\varphi_{1}}{\longrightarrow}Y_{2}\overset{\varphi_{2}}%
{\longrightarrow}\cdots\overset{\varphi_{n-1}}{\longrightarrow}Y_{n}%
\overset{\varphi_{n}}{\longrightarrow}Y_{1}$, $Y_{i}=\mathfrak{F}%
_{\operatorname*{Ob}}^{(n)}(X_{i})\in\mathcal{D}_{\operatorname*{Ob}%
},\,\varphi_{i}=\mathfrak{F}_{\operatorname*{Mor}}^{(n)}(X_{i})\in
\mathcal{D}_{\operatorname*{Mor}}$.
\end{lemma}

\textbf{Proof:} Because $\mathfrak{F}_{\operatorname*{Ob}}^{(n)}%
=\mathfrak{F}_{\operatorname*{Ob}}$, we need to prove that $\mathfrak{F}%
_{\operatorname*{Mor}}^{(n)}(\psi_{i})\circ\mathbf{e}_{X_{i}}^{(n)}%
=\mathfrak{F}_{\operatorname*{Mor}}^{(n)}(\psi_{i})$. Indeed, we have
\[
\mathfrak{F}_{\operatorname*{Mor}}^{(n)}(\psi_{i})=\mathfrak{F}%
_{\operatorname*{Mor}}^{(n)}\left(  \psi_{i}\circ\mathbf{e}_{X_{i}}%
^{(n)}\right)  =\mathfrak{F}_{\operatorname*{Mor}}^{(n)}\left(  \psi
_{i}\right)  \circ\mathfrak{F}_{\operatorname*{Mor}}^{(n)}\left(  e_{X_{i}%
}^{(n)}\right)  =\mathfrak{F}_{\operatorname*{Mor}}^{(n)}(\psi_{i})\circ
e_{X_{i}}^{(n)}.
\]
\hfill$\Box$\newline 

Multifunctors can be regularized in a similar way. Next we ``regularize''
natural transformations. Let $\mathfrak{F}^{(n)}$ and $\mathfrak{G}^{(n)}$ be
two $n$-regular functors of the category $\mathcal{C}$ into the category
$\mathcal{D}$.

\begin{definition}
A natural transformation $\mathfrak{s}:\mathfrak{F}^{(n)}\rightarrow
\mathfrak{G}^{(n)}$ is a collection of morphisms $\mathfrak{s}=\{\mathfrak{s}%
_{X}:\mathfrak{F}_{\operatorname*{Ob}}(X)\rightarrow G_{\operatorname*{Ob}%
}(X),\;X\in\mathcal{C}_{\operatorname*{Ob}}\}$ such that
\begin{equation}
\mathfrak{s}_{Y}\circ\mathfrak{F}_{\operatorname*{Mor}}^{(n)}(\alpha
)=\mathfrak{G}_{\operatorname*{Mor}}^{(n)}(\alpha)\circ\mathfrak{s}_{X},
\end{equation}
for every regular morphism $\alpha:X\rightarrow Y$,is said to be
$n$-\textit{regular natural transformation}.
\end{definition}

\begin{definition}
A monoidal category $\mathcal{C}\equiv\mathcal{C}(\otimes,\mathbb{K})$
equipped with a family of obstruction morphisms $\mathbf{E}^{(n)}%
=\{\mathbf{e}_{X}^{\left(  n\right)  }:X\in\mathcal{C}_{\operatorname*{Ob}%
};n=1,2,...\}$ satisfying the ``homomorphism'' condition
\begin{equation}
\mathbf{e}_{X\otimes Y}^{(n)}=\mathbf{e}_{X}^{(n)}\otimes\mathbf{e}_{Y}%
^{(n)}.\label{mul}%
\end{equation}
is said to be an \textit{obstructed monoidal category}.
\end{definition}

Let $\mathcal{C}$ be a monoidal category with duality \cite{yet,maj3}. If
$(X,\psi)$ a regular $n$-cocycle on $\mathcal{C}$, then there is a
corresponding regular $n$-cocycle $(X^{\vee},\psi^{\vee})$, called the dual of
$(X,\psi)$. This means that we have a sequence
\begin{equation}%
\begin{array}
[c]{c}%
X_{1}^{\vee}\overset{\psi_{n}^{\vee}}{\rightarrow}X_{n}^{\vee}\overset
{\psi_{n-1}^{\vee}}{\rightarrow}\cdots\overset{\psi_{2}^{\vee}}{\rightarrow
}X_{2}^{\vee}\overset{\psi_{1}^{\vee}}{\rightarrow}X_{1}^{\vee}%
\end{array}
\label{regub}%
\end{equation}
such that there is an (endo-)morphism $\mathbf{e}_{X_{1}^{\vee}}^{\left(
n\right)  }:X_{1}^{\vee}\rightarrow X_{1}^{\vee}$ defined by the following
equation
\begin{equation}%
\begin{array}
[c]{c}%
\mathbf{e}_{X_{1}^{\vee}}^{\left(  n\right)  }:=\psi_{n}^{\vee}\circ
\cdots\circ\psi_{2}^{\vee}\circ\psi_{1}^{\vee}%
\end{array}
\label{egub}%
\end{equation}
and subjects to the following relation
\begin{equation}
\psi_{1}^{\vee}\circ\psi_{n}^{\vee}\circ\cdots\circ\psi_{2}^{\vee}\circ
\psi_{1}^{\vee}=\psi_{1}^{\vee},
\end{equation}
where $X_{i}^{\vee}=\mathcal{C}(X_{i},\mathbb{K})$, $\psi_{i}^{\vee}%
:X_{i+1}^{\vee}\rightarrow X_{i}^{\vee}$, $i=1,\ldots,n-1$. Let $\mathsf{g}%
_{X}\equiv\langle\ldots|\ldots\rangle_{X}:X^{\vee}\otimes X\rightarrow
\mathbb{K}$ be a pairing in $\mathcal{C}$.

\begin{definition}
If for every object $X\in\mathcal{C}_{\operatorname*{Ob}}$ with an obstruction
$\mathbf{e}_{X}^{(n)}$ there is is an object $X^{\vee}$ with an obstruction
$\mathbf{e}_{X^{\vee}}^{(n)}$ such that
\begin{equation}%
\begin{array}
[c]{l}%
\langle\mathbf{e}_{X^{\vee}}^{(n)}(x^{\vee}),x\rangle=\langle x^{\vee
}|\mathbf{e}_{X}^{(n)}(x)\rangle_{\mathcal{U}},
\end{array}
\label{duality}%
\end{equation}
for $x^{\vee}\in X^{\vee},x\in X$, then the category $\mathcal{C}$ is said to
be an \textit{obstructed category with duality}.
\end{definition}

The category $\mathcal{C}_{\mathcal{A}}^{\left(  2\right)  }$ above is an
obstructed category with duality, where.
\begin{equation}%
\begin{array}
[c]{cc}%
X_{1}^{\vee}:=\mathcal{C}_{\mathcal{A}}^{\left(  2\right)  }(X_{1}%
,\mathbb{K}), & X_{2}^{\vee}:=\mathcal{C}_{\mathcal{A}}^{\left(  2\right)
}(X_{2},\mathbb{K}),
\end{array}
\end{equation}
are spaces of linear functionals on $X_{1}$ and $X_{2}$, respectively, and
\begin{equation}%
\begin{array}
[c]{cc}%
X_{1}^{\vee}:=lin_{\mathbb{K}}\{\Xi_{1},\Xi_{2}\Xi_{1}\}, & X_{2}^{\vee
}:=lin_{\mathbb{K}}\{\Xi_{2},\Xi_{1}\Xi_{2}\}.
\end{array}
\end{equation}
where
\begin{equation}%
\begin{array}
[c]{cc}%
\langle\Xi_{1}|\Theta_{1}\rangle_{X_{1}}:=1, & \langle\Xi_{2}|\Theta
_{2}\rangle_{X_{2}}:=1
\end{array}
\end{equation}
and
\begin{equation}%
\begin{array}
[c]{cc}%
\langle\Xi_{1}\Xi_{2}|\Theta_{2}\Theta_{1}\rangle_{X_{1}}:=1 & \langle\Xi
_{2}\Xi_{1}|\Theta_{1}\Theta_{2}\rangle_{X_{2}}:=1.
\end{array}
\end{equation}
The monoidal product is the $\mathbb{K}$-tensor product \cite{yet,maj3}.

\section{Regular algebras and bialgebras}

Let $\mathcal{C}$ be an obstructed monoidal category \cite{dup/mar3,dup/mar5}.

\begin{definition}
An algebra $\mathcal{A}^{reg}$ in the category $\mathcal{C}$ such that the
multiplication $\mathrm{m}:\mathcal{A}^{reg}\otimes\mathcal{A}^{reg}%
\rightarrow\mathcal{A}^{reg}$ is a regular $n$-cocycle morphism
\begin{equation}
\mathrm{m}\circ(\mathbf{e}_{\mathcal{A}}^{(n)}\otimes\mathbf{e}_{\mathcal{A}%
}^{(n)})=\mathbf{e}_{\mathcal{A}}^{(n)}\circ\mathrm{m},\label{regal}%
\end{equation}
is said to be a \textit{regular }$n$\textit{-cocycle algebra}.
\end{definition}

\begin{lemma}
The graded algebra $\mathcal{A}^{reg}\left(  \Theta_{1},\Theta_{2}\right)  $
is a regular $2$-cocycle algebra in the category $\mathcal{C}_{\mathcal{A}}$.
\end{lemma}

\textbf{Proof:} Let $\mathbf{a},\mathbf{b}\in\mathcal{A}^{reg}\left(
\Theta_{1},\Theta_{2}\right)  $ have the form (\ref{a}). Then for the
obstruction $\mathbf{e}_{\mathcal{A}}^{(2)}$ we calculate
\begin{equation}
\mathbf{e}_{\mathcal{A}}^{(2)}(\mathbf{a}):=1+a_{2}\Theta_{1}+a_{1}\Theta
_{2}+a_{21}\Theta_{1}\Theta_{2}+a_{12}\Theta_{2}\Theta_{1} \label{e2}%
\end{equation}
So the condition (\ref{regal}) for $\mathbf{e}_{\mathcal{A}}^{(2)}%
(\mathbf{a}),\mathbf{e}_{\mathcal{A}}^{(2)}(\mathbf{b}),\mathbf{e}%
_{\mathcal{A}}^{(2)}(\mathbf{c})$ in components is%
\begin{align*}
c_{1}  &  =a_{1}+b_{1}+a_{2}b_{12}+a_{21}b_{2},\\
c_{2}  &  =a_{2}+b_{2}+a_{1}b_{21}+a_{12}b_{1},\\
c_{12}  &  =a_{12}+b_{12}+a_{1}b_{2}+a_{12}b_{12},\\
c_{21}  &  =a_{21}+b_{21}+a_{2}b_{1}+a_{21}b_{21}.
\end{align*}
. \hfill$\Box$

Note that in case of ground field is complex among obstructions $\mathbf{e}%
_{\mathcal{A}}^{(2)}(\mathbf{a})$ there are exist only two ones which are
idempotent as elements, i. e. satisfy the equation $\mathbf{e}_{\mathcal{A}%
}^{(2)}(\mathbf{a})\mathbf{e}_{\mathcal{A}}^{(2)}(\mathbf{a})=\mathbf{e}%
_{\mathcal{A}}^{(2)}(\mathbf{a})$, and they do not depend from $\mathbf{a}$%
\[
\mathbf{e}_{\mathcal{A}}^{\left(  2\right)  \text{idemp}}=1+\Theta_{1}%
+\Theta_{2}+\left(  -\dfrac{1}{2}\pm i\dfrac{\sqrt{3}}{2}\right)  \Theta
_{1}\Theta_{2}+\left(  -\dfrac{1}{2}\mp i\dfrac{\sqrt{3}}{2}\right)
\Theta_{2}\Theta_{1}.
\]

There is a graded algebra $\mathcal{A}^{reg,\vee}(\Xi_{1},\Xi_{2})$ generated
by $\Xi_{1}$ and $\Xi_{2}$ subject to the relations
\begin{align}
\Xi_{1}^{2}  &  =\Xi_{2}^{2}=0,\nonumber\\
\Xi_{1}\Xi_{2}\Xi_{2}  &  =\Xi_{1},\;\;\;\Xi_{2}\Xi_{1}\Xi_{2}=\Xi_{2}.
\label{drel}%
\end{align}
which is also is a regular $2$-cocycle algebra in the category $\mathcal{C}$.

\hfill$\Box$

\begin{definition}
A coalgebra $\mathcal{A}^{reg}$ in the category $\mathcal{C}$ such that the
comultiplication $\triangle_{reg}:\mathcal{A}^{reg}\longrightarrow
\mathcal{A}^{reg}\otimes\mathcal{A}^{reg}$ is a regular $n$-cocycle morphism
\begin{equation}
\triangle_{reg}\circ\mathbf{e}_{\mathcal{A}}^{(n)}=(\mathbf{e}_{\mathcal{A}%
}^{(n)}\otimes\mathbf{e}_{\mathcal{A}}^{(n)})\circ\triangle_{reg}%
.\label{regco}%
\end{equation}
is said to be a \textit{regular }$n$\textit{-cocycle coalgebra}.
\end{definition}

\begin{lemma}
Let $\mathcal{A}$ be an algebra in an obstructed category with duality, and
let $\mathcal{A^{\vee}}$ be its dual coalgebra%
\begin{equation}%
\begin{array}
[c]{l}%
\langle\triangle_{reg}(\mathbf{a}^{\vee}),\mathbf{a}_{1}\otimes\mathbf{a}%
_{2}\rangle=\langle\mathbf{a}^{\vee},\mathrm{m}(\mathbf{a}_{1}\otimes
\mathbf{a}_{2})\rangle,
\end{array}
\label{dual3}%
\end{equation}
where $\mathbf{a}_{1},\mathbf{a}_{2}\in\mathcal{A},\;\mathbf{a}^{\vee}%
\in\mathcal{A}^{\vee}$. If $\mathcal{A}$ is a regular $n$-cocycle algebra,
then $\mathcal{A}^{\vee}$ is a regular n-cocycle coalgebra.
\end{lemma}

\textbf{Proof:} Let us apply the regularity condition (\ref{regal}) to the
above duality condition (\ref{dual3}). Then the lemma follows from relations
(\ref{mul}),(\ref{duality}), and (\ref{regco}). \hfill$\Box$

The graded algebra $\mathcal{A}^{\vee}(\Xi_{1},\Xi_{2})$ is also a regular
$2$-cocycle coalgebra which is dual to $\mathcal{A}^{reg}(\Theta_{1}$
$,\Theta_{2}$ $)$ in the category $\mathcal{C}$.

\begin{definition}
Let $\mathcal{A}^{reg}$ be a regular $n$-cocycle algebra. If $\mathcal{A}%
^{reg}$ is also regular coalgebra such that $\Delta_{reg}\left(
\mathbf{ab}\right)  =\Delta_{reg}\left(  \mathbf{a}\right)  \Delta
_{reg}\left(  \mathbf{b}\right)  $, then it is said to be a \textit{regular
}$n$\textit{-cocycle almost bialgebra}.
\end{definition}

One can also consider the corresponding generalization of Hopf algebras
\cite{nil1,fangli3}.

The graded algebra $\mathcal{A}^{reg}(\Theta_{1}$ $,\Theta_{2}$ $)$ is a
regular $2$-cocycle almost bialgebra, where
\begin{align}
\triangle_{reg}(\Theta_{1})  &  =\Theta_{1}\otimes\mathbf{e}_{\mathtt{X}_{1}%
}+\mathbf{e}_{\mathtt{X}_{1}}\otimes\Theta_{1},\\
\triangle_{reg}(\Theta_{2})  &  =\Theta_{2}\otimes\mathbf{e}_{\mathtt{X}_{2}%
}+\mathbf{e}_{\mathtt{X}_{2}}\otimes\Theta_{2},
\end{align}

Similarly, the graded algebra $\mathcal{A}^{\vee}(\Xi_{1},\Xi_{2})$ is a
regular $2$-cocycle almost bialgebra.

\begin{definition}
Let $M^{reg}$ be an $\mathcal{A}^{reg}$-module equipped with the
$\mathcal{A}^{reg}$-module action $\rho_{M}^{reg}:\mathcal{A}^{reg}\otimes
M^{reg}\longrightarrow M^{reg}$, where $\mathcal{A}^{reg}$ is a $n$-regular
cocycle algebra. Then $M^{reg}$ is said to be an $n$-regular cocycle
$\mathcal{A}^{reg}$-module, if an only if
\begin{equation}
\rho_{M}^{reg}\circ(\mathbf{e}_{\mathcal{A}}^{(n)}\otimes\mathbf{e}_{M}%
^{(n)})=\mathbf{e}_{M}^{(n)}\circ\rho_{M}^{reg}. \label{rmod}%
\end{equation}
\end{definition}

\section{Regular Wick algebras}

A pair of algebras $\mathcal{A}$ and $\mathcal{A}^{\dag}$ equipped with an
antilinear and involutive anti-isomorphism $(-)^{\dag}:\mathcal{A}%
\rightarrow\mathcal{A}^{\dag}$ such that
\begin{equation}
\mathrm{m}_{\mathcal{A}^{\dag}}(\mathbf{b}^{\dag}\otimes\mathbf{a}^{\dag
})=(\mathrm{m}_{\mathcal{A}}(\mathbf{a}\otimes\mathbf{b}))^{\dag}%
,\quad(\mathbf{a}^{\dag})^{\dag}=\mathbf{a},
\end{equation}
where $\mathbf{a},\mathbf{b}\in\mathcal{A}$ and $\mathbf{a}^{\dag}%
,\mathbf{b}^{\dag}\in\mathcal{A}^{\dag}$ are their images under the
isomorphism $(-)^{\dag}$, is said to be the \textit{conjugated} one to the
other \cite{bor/mar}.

\begin{definition}
Let $(\mathcal{A},\mathcal{A}^{\dag})$ be a pair of conjugated algebras. A
linear and mapping $\Psi:\mathcal{A}^{\dag}\otimes\mathcal{A}\rightarrow
\mathcal{A}\otimes\mathcal{A}^{\dag}$ such that
\begin{equation}%
\begin{array}
[c]{l}%
\Psi\circ(\operatorname*{id}\nolimits_{\mathcal{A}^{\dag}}\otimes
\mathrm{m}_{\mathcal{A}})=(\mathrm{m}_{\mathcal{A}}\otimes\operatorname*{id}%
\nolimits_{\mathcal{A}^{\dag}})\circ(\operatorname*{id}\nolimits_{\mathcal{A}%
}\otimes\Psi)\circ(\Psi\otimes\operatorname*{id}\nolimits_{\mathcal{A}}),\\
\Psi\circ(\mathrm{m}_{\mathcal{A}^{\dag}}\otimes\operatorname*{id}%
\nolimits_{\mathcal{A}})=(\operatorname*{id}\nolimits_{\mathcal{A}}%
\otimes\mathrm{m}_{\mathcal{A}^{\dag}})\circ(\Psi\otimes\operatorname*{id}%
\nolimits_{\mathcal{A}^{\dag}})\circ(\operatorname*{id}\nolimits_{\mathcal{A}%
^{\dag}}\otimes\Psi)
\end{array}
\label{twc}%
\end{equation}
is said to be a \textit{cross symmetry} or \textit{generalized twist}
\cite{bor/mar}.
\end{definition}

\begin{definition}
The tensor product $\mathcal{A}\otimes\mathcal{A}^{\dag}$ equipped with the
multiplication
\begin{equation}%
\begin{array}
[c]{c}%
\mathrm{m}_{\Psi}:=(\mathrm{m}_{\mathcal{A}}\otimes\mathrm{m}_{\mathcal{A}%
^{\dag}})\circ(\operatorname*{id}\nolimits_{\mathcal{A}}\otimes\Psi
\otimes\operatorname*{id}\nolimits_{\mathcal{A}^{\dag}})
\end{array}
\label{mult}%
\end{equation}
is an associative algebra called a \textit{Hermitian Wick algebra}
\cite{bor/mar,jor/sch/wer} and it is denoted by $\mathcal{W}=\mathcal{W}%
_{\Psi}(\mathcal{A})=\mathcal{A}>\!\!\!\vartriangleleft_{\Psi}\mathcal{A}^{\dag}$.
\end{definition}

\begin{definition}
Let $(\mathcal{A},\mathcal{A}^{\dag})$ be a pair of conjugated algebras, and
let $\Psi:\mathcal{A}^{\dag}\otimes\mathcal{A}\rightarrow\mathcal{A}%
\otimes\mathcal{A}^{\dag}$ be a cross symmetry. If both algebras $\mathcal{A}$
and $\mathcal{A}^{\dag}$ are $n$-regular cocycle algebras and
\begin{equation}%
\begin{array}
[c]{c}%
(\mathbf{e}_{\mathcal{A}}^{(n)}\otimes\mathbf{e}_{\mathcal{A}^{\dag}}%
^{(n)})\circ\Psi^{reg}:=\Psi^{reg}\circ(\mathbf{e}_{\mathcal{A}^{\dag}}%
^{(n)}\otimes\mathbf{e}_{\mathcal{A}}^{(n)}),
\end{array}
\end{equation}
then the mapping $\Psi^{reg}$ is said to be $n$-\textit{regular cocycle cross
symmetry} and the corresponding algebra $\mathcal{W}^{reg}=\mathcal{W}%
_{\Psi^{reg}}(\mathcal{A})=\mathcal{A}>\!\!\!\vartriangleleft_{\Psi}\mathcal{A}^{\dag}$ is
called $n$-\textit{regular cocycle Wick algebra}.
\end{definition}

A natural regular generalization of (\ref{mult}) is given by

\begin{definition}
The tensor product $\mathcal{A}\otimes\mathcal{A}^{\dag}=\mathcal{A}%
^{reg}\otimes\mathcal{A}^{reg,\dag}$ equipped with the multiplication
\begin{equation}%
\begin{array}
[c]{c}%
\mathrm{m}_{\Psi}^{reg}:=(\mathrm{m}_{\mathcal{A}}\otimes\mathrm{m}%
_{\mathcal{A}^{\dag}})\circ(\mathbf{e}_{\mathcal{A}}^{(n)}\otimes\Psi
^{reg}\otimes\mathbf{e}_{\mathcal{A}^{\dag}}^{(n)})
\end{array}
\end{equation}
is an associative algebra called a \textit{regular Hermitian Wick algebra} and
it is denoted by $\mathcal{W}^{reg}=\mathcal{W}_{\Psi^{reg}}(\mathcal{A}%
)=\mathcal{A}>\!\!\!\vartriangleleft_{\Psi}\mathcal{A}^{\dag}$.
\end{definition}

Let us describe the regular Wick algebra corresponding to the algebra
$\mathcal{A}^{reg}(\Theta_{1},\Theta_{2})$. In this case $\mathcal{A}^{\dag
}=\mathcal{A}^{reg,\vee}(\Xi_{1},\Xi_{2})$. Let us assume that
\begin{equation}%
\begin{array}
[c]{ll}%
\Psi^{reg}(\Xi_{1}\otimes\Theta_{1})=e_{\mathtt{X}_{1}}-\Theta_{1}\otimes
\Xi_{1} & \Psi^{reg}(\Xi_{1}\otimes\Theta_{2})=\Theta_{2}\otimes\Xi_{1},\\
\Psi^{reg}(\Xi_{2}\otimes\Theta_{2})=e_{\mathtt{X}_{2}}-\Theta_{2}\otimes
\Xi_{2} & \Psi^{reg}(\Xi_{2}\otimes\Theta_{1})=\Theta_{1}\otimes\Xi_{2},
\end{array}
\label{p1}%
\end{equation}
where $e_{\mathtt{X}_{i}}$ are defined in (\ref{eff}). We use the (\ref{twc})
relations in order to calculate $\Psi^{reg}(\Xi_{1}\otimes\Theta_{1}$
$\Theta_{2}$ $)$ and others similar expressions. We obtain
\begin{equation}%
\begin{array}
[c]{l}%
\Psi^{reg}(\Xi_{1}\otimes\Theta_{1}\Theta_{2})=\Psi^{reg}(\Psi^{reg}(\Xi
_{1}\otimes\Theta_{1})\otimes\Theta_{2})\\
=\Theta_{2}-\Theta_{1}\Psi^{reg}(\Xi_{1}\otimes\Theta_{2})=\Theta_{2}%
-\Theta_{1}\Theta_{2}\otimes\Xi_{1}%
\end{array}
\label{p2}%
\end{equation}

Thus using (\ref{p1}) and (\ref{p2}) we can construct an $n$-regular cocycle
Wick algebra $\mathcal{W}^{reg}=\mathcal{W}_{\Psi^{reg}}\left(  \mathcal{A}%
^{reg}\left(  \Theta_{1},\Theta_{2}\right)  \right)  $.

\textbf{Acknowledgements}. S.D. is grateful to Alexander Stolin and Leonid
Vaksman for useful discussions and to Jerzy Lukierski for kind hospitality at
the University of Wroclaw, where this work was finished. The work of W.M. is
partially sponsored by Polish Committee for Scientific Research (KBN) under
Grant No 5P03B05620.

\end{document}